\documentclass[11pt]{amsart}

\usepackage{hyperref}

 \newcommand{\doubleheaddownarrow}{\big\downarrow\kern-3.325mm\downarrow}

 \newcommand{\iso}{\cong}

\newcommand{\R}{\mathbb R}
 \newcommand{\Z}{\mathbb Z}

 \newcommand{\D}{\Delta}

 \newcommand{\Tor}{\operatorname{Tor}}

 \newcommand{\rank}{\operatorname{rank}}

 \newcommand{\pf}{\operatorname{pf}}

 \newcommand{\Hom}{\operatorname{Hom}}

 \newcommand{\Spec}{\operatorname{Spec}}

\newcommand{\romstar}{\operatorname{star}}

\newtheorem{theorem}{Theorem}[section]
\newtheorem{lemma}[theorem]{Lemma}
\newtheorem{prop}[theorem]{Proposition}

\theoremstyle{definition}

\newtheorem{example}[theorem]{Example}

\theoremstyle{remark}
\newtheorem{remark}[theorem]{Remark}

\numberwithin{equation}{section}

\begin{document}

\title [Stanley-Reisner rings associated to  cyclic polytopes]{On the structure of 
                    Stanley-Reisner rings associated to cyclic polytopes}

\author{Janko B\"{o}hm}
\address{Department of Mathematics, University of California \\ Berkeley
\\ CA 94720, USA, and Department of Mathematics, Universit\"{a}t des
Saarlandes, Campus E2 4 \\ D-66123 \\ Saarbr\"{u}cken, Germany}
\email{boehm@math.uni-sb.de}

\thanks{J. B. supported by DFG (German Research Foundation) through Grant BO3330/1-1.
S. P. was supported by the Portuguese Funda\c{c}\~ao para a Ci\^encia e a Tecno\-lo\-gia through 
Grant SFRH/BPD/22846/2005 of POCI2010/FEDER and  through Project PTDC/MAT/099275/2008.
}

\author{Stavros Argyrios Papadakis}
\address{Centro de An\'{a}lise Matem\'{a}tica, Geometria e Sistemas Din\^{a}micos,
Departamento de Matem\'atica, Instituto Superior T\'ecnico, 
Universidade T\'ecnica de Lisboa,
Av. Rovisco Pais, 1049-001 Lisboa,
Portugal}
\email{papadak@math.ist.utl.pt}


\begin{abstract} 

We study the structure of Stanley--Reisner rings 
associated to cyclic polytopes, using ideas from unprojection theory.
Consider the boundary simplicial complex  $\Delta (d,m)$  of
the $d$-dimensional cyclic polytope with $m$ vertices.
We show how to express the  Stanley-Reisner ring of $\Delta (d,m+1)$  
in terms of the Stanley--Reisner rings of  $\Delta (d,m)$ and 
$\Delta (d-2,m-1)$.
 As an application, we use the Kustin--Miller complex construction to
identify the minimal graded free resolutions of these rings. In particular,
we recover results of Schenzel, Terai and Hibi about their graded Betti numbers.
\end{abstract}

\maketitle

\section { Introduction }

Gorenstein commutative rings form an important class of commutative rings.
For example, they appear in algebraic geometry as canonical rings 
of regular surfaces and  anticanonical rings of  Fano $n$-folds and
in algebraic combinatorics as Stanley--Reisner rings of sphere
triangulations.  
In codimensions 
$1$ and $2$ they are complete intersections and in codimension $3$ they
are Pfaffians \cite{BE}, but, to our knowledge, no structure theorems
are known for higher codimensions.

Unprojection theory  \cite{R1},   which analyzes 
and constructs complicated commutative rings in terms of simpler ones, 
began with the aim of  partly filling this gap.  The first kind of unprojection
which appeared in the literature is that of type Kustin--Miller,
studied originally by Kustin and Miller \cite{KM} and later by
Reid and the second author \cite{P1, PR}. Starting from a
codimension $1$ ideal $J$ of a Gorenstein ring $R$ such that the
quotient $R/J$ is Gorenstein, Kustin--Miller unprojection uses the
information contained in $\Hom_{R}(J,R)$ to construct a new
Gorenstein ring $S$ which is birational to $R$ and corresponds
to the contraction of $V(J) \subset \Spec R$. See
Subsection~\ref{subs!defnofkustinmillerunprojection} for a precise 
definition of Kustin--Miller unprojection and the introduction
of \cite{BP} for references to applications.

In the paper \cite{BP}, the authors proved that on the algebraic 
level of Stanley--Reisner rings, stellar subdivisions of Gorenstein* 
simplicial complexes correspond to  Kustin--Miller unprojections
and gave applications to Stanley-Reisner rings associated to 
stacked polytopes.  In the present paper, we use unprojection theory to  
study the  structure of Stanley--Reisner rings associated to cyclic 
polytopes. This setting is different from the one studied in \cite{BP} 
since here, except for some easy subcases, stellar subdivisions do not appear
and the unprojection ideals are more complicated.

Our main result, which is stated precisely in 
Theorems~\ref{theorem!mainthm_cyclicpolytopes_deven} and 
\ref{theorem!mainthm_cyclicpolytopes_dodd}, can be described  as follows.
Assume $d \geq 4$ and $d+1 < m$.  Consider the  cyclic 
polytope which has $m$ vertices  and   dimension $d$, 
and   denote by $\Delta (d,m)$ its
boundary simplicial complex. We show how to express the Stanley-Reisner ring 
of $\Delta (d,m+1)$ 
in terms of the Stanley--Reisner rings of $\Delta (d,m)$ and $\Delta (d-2,m-1)$
via  Kustin--Miller unprojection.
Moreover, a similar result is also true for the remaining cases 
$d=2,3$ and $m=d+1$, see Subsections~\ref{sub!thecasedis2},
\ref{sub!casedevenandmisdplus}, \ref{sub!casedis3and} and
\ref{sub!casedoddandmisdplusone}. In Section~\ref{sec!combinatorial_interpretation}
we give a combinatorial interpretation of our construction.

As an application, in Section~\ref{sec!application_to_resolutions}
we inductively identify the minimal graded free resolutions of the
Stanley--Reisner rings $k[\Delta (d,m)]$.
We use this identification in Proposition~\ref {prop!keytechnicalforapplication}
to calculate the graded Betti numbers of these rings,  recovering 
results originally due to Schenzel \cite{S} for $d$ even and  Terai and Hibi \cite{TH} 
for $d$ odd. Our derivation is
more algebraic than the one in \cite{TH}, and does not use Hochster's formula 
or Alexander duality.  Finally, Subsection~\ref{subs!examplesforKMconstruction}  
contains examples and a link to related computer algebra code.

An interesting open question is whether there
are other families of Gorenstein Stanley--Reisner rings related by
unprojections in a similar way as cyclic polytopes, compare also the
discussion in \cite[Section~6]{BP}.

\section { Preliminaries  }  \label{sec!preliminaries}

Assume $k$ is a field, and $m$  a positive integer. An
(abstract) simplicial complex on the vertex set $\{1, \dots ,m \}$ is a
collection $\D$ of subsets of $\{1, \dots ,m \}$ such that (i) all singletons
$\{i\}$ with $i \in \{1, \dots ,m\}$ belong to $\D$ and (ii) $\sigma \subset
\tau \in \D$ implies $\sigma \in \D$.  The elements of $\D$ are
called \emph{faces} and those maximal with respect to inclusion
are called \emph{facets}. The dimension of a face $\sigma$ is
defined as one less than the cardinality of $\sigma$. The
dimension of $\D$ is the maximum dimension of a face. 
Any abstract simplicial complex $\Delta$ has a
geometric realization, which is unique up to linear homeomorphism.

For any subset $W$ of $\{1, \dots , m \}$, we denote by $x_W$ the
square-free monomial in the polynomial ring $k[x_1,\dots,x_m]$
with support $W$, in other words $x_W$ is the product of $x_t$ for  $t \in W$.
The ideal $I_\D$ of $k[x_1,\dots,x_m]$ which is generated by the square-free 
monomials $x_W$ with $W\notin \D$ is called the \emph{Stanley-Reisner ideal} of $\D$. The
\emph{face ring}, or \emph{Stanley-Reisner ring}, of $\D$ over $k$, 
denoted $k[\D]$, is defined as the quotient ring of $k[x_1,\dots,x_m]$ by the ideal
$I_\D$.

Assume $R=k[x_1, \dots ,x_m]$ is a polynomial ring over a field $k$
with the degrees of all variables $x_i$ positive,   and denote by
$m=(x_1, \dots ,x_m)$ the maximal homogeneous ideal of $R$. Assume
$M$ is a finitely  generated graded $R$-module. Denote by
\[
      0 \to F_{g} \to F_{g-1} \to \dots  \to F_1 \to F_0 \to M\to 0
\] 
the minimal graded free resolution of $M$ as $R$-module, and write  
\[   
         F_i = \oplus_{j} R(-j)^{b_{ij}}.
\]
The integer  $b_{ij}$  is called the \emph{$ij$-th graded Betti number} of $M$ 
and is also  denoted by $b_{ij}(M)$.
For fixed $i$ we set $b_i(M) = \sum_{j} b_{ij}(M)$. The integer $b_i(M)$ is 
the rank of the free $R$-module $F_i$ in the category of (ungraded) $R$-modules, and
\begin {equation} \label{eqn!total_betti_as_tor_dimensions}
   b_{i}(M) = \dim_{R/m} \Tor_i^{R}(R/m, M),
\end{equation} cf. \cite[Proposition~1.7]{Ei2}.
For more details  about free resolutions and Betti numbers 
see, for example, \cite[Sections~19, 20] {Ei}.

Assume $R$ is a ring. An element $r \in R$ will be called \emph{$R$-regular} if
the multiplication by $r$ map $R \to R, u \mapsto ru$ is injective. A sequence 
$r_1, \dots ,r_n$ of elements of
$R$ will be called a \emph{regular $R$-sequence} if $r_1$ is $R$-regular,   and,
for $2 \leq i \leq n$, we have that  $r_i$ is $R/(r_1, \dots ,r_{i-1})$-regular.

Assume $k$ is a field, and    $a,m,n$ three positive  integers   with   $m <  n$ and   $2a  \leq  n-m+2$.
  We define the ideal   $I_{a,m,n}   \subset  k[ x_m, x_{m+1}, \dots ,  x_n]$ by 
  \[
        I_{a,m,n} = (  x_{t_1}x_{t_2} \dots  x_{t_a}  \bigm| m \leq t_1, t_a \leq n, 
         t_{j}+2 \leq t_{j+1}  \text{ for } 1 \leq j \leq a-1 ).
  \]
 The assumption   $2a  \leq  n-m+2$ implies that
 there exists at least one monomial generator of $I_{a,m,n}$, namely $x_mx_{m+2} \dots  x_{m+2(a-1)}$.
 For example, we have  $I_{2,3,6} = (x_3x_5,x_3x_6,x_4x_6)$.

\subsection {Cyclic polytopes}   \label{subs!cyclic_polytopes_dfn}

   Recall from \cite [Section~5.2]{BH} the definition of cyclic polytopes.
   We fix two integers $m,d$, with  $2 \leq d <  m$, and  define the cyclic polytope
             $C_d (m)   \subset \R^d$    
  as follows:  Fix, for $1 \leq i \leq m$,   $t_i \in \R$   with   $t_1 < t_2  < \dots  < t_m$. By definition,
  the cyclic polytope  $C_d (m)  =  C_d (t_1, \dots, t_m)$ is the convex hull in  $\R^d$ of the subset
 $\{  f(t_1), f(t_2),  \dots , f (t_m)    \} \subset \R^d$,
  where   $f \colon \R \to \R^d$ with $f (t) = (t, t^2, \dots , t^d)$  for   $t \in \R$.
  We have that $C_d(m)$ is a simplicial $d$-polytope, which 
  up to combinatorial equivalence does not depend on the choice of the points $t_i$.
  We denote by  $\Delta (d,m)$ the  boundary simplicial complex of $C_d (m)$,
  by definition  $\Delta (d,m)$ has as  elements the empty set and the sets of  vertices of the proper faces of $C_d(m)$, 
  cf.~\cite[Corollary~5.2.7]{BH}.

  Assume $W \subset \{1, \dots ,m \}$ is a proper nonempty subset. A nonempty subset $X \subset W$ is called
  contiguous  if there exist $i,j$ with $2 \leq i \leq j \leq m-1$ such that $i-1 \notin W$, $j+1 \notin W$,
  $X = \{i,i+1, \dots , j \}$.    A contiguous $X \subset W$ is called odd contiguous if $\# X$ is odd.
  Assume $W$ contains a contiguous subset, this is equivalent to the existence of 
  $a \in W$ and $b_1, b_2 \in \{1, \dots ,m \} \setminus W$ with $b_1 < a < b_2$.
  Then, there exist a unique integer $t \geq 1$ and
  a unique decomposition 
  \[
     W  = Y_1 \cup X_1 \cup X_2 \cup  \dots  \cup X_t \cup Y_2,
  \]
  such that $Y_1$ is either empty or of the form $\{1, 2, \dots , i \}$  for some
  $i \geq 1$ with $i+1 \notin W$,    $Y_2$ is either empty or of the form $\{j, j+1, \dots , m \}$
  for some $j \leq m$ with $j-1 \notin W$, each $X_p$, for $1 \leq p \leq t$, is a contiguous subset of
  $W$, and for  $p_1 < p_2$ each element of $X_{p_1}$ is strictly smaller than any element of $X_{p_2}$.

  For a real number $r$ we denote by $[r]$ the integral value of $r$,
  i.e., the largest integer which is smaller or equal than $r$.
  The following theorem characterizing the faces of $\Delta (d,m)$ is proven in \cite[Theorem~5.2.13]{BH}, 
  compare also \cite[Lemma~2.2]{TH}. 
  
 \begin {theorem}  \label {thm!facesofcyclicpolytopinfo}
  Assume $W \subset \{ 1, \dots ,m \}$ is a  nonempty subset with $\# W \leq d$.
  $W$ is a face of $\Delta (d,m)$ if and only if the number of odd contiguous subsets of $W$ is at most $d-\# W$.
  In particular, if $\# W  \leq [d/2]$ then  $W$ is a face of $\Delta (d,m)$.  
\end {theorem}

\subsection {Kustin--Miller unprojection} \label{subs!defnofkustinmillerunprojection}

We recall the definition of Kustin--Miller unprojection from \cite{PR}. 
Assume $R$ is a local (or graded) Gorenstein ring, and $J \subset R$ a codimension $1$
ideal with $R/J$ Gorenstein. Fix $\phi \in \Hom_R(J,R)$ such that $\Hom_R(J,R)$  is
generated as an $R$-module by the subset $\{i, \phi \}$, where $i$ denotes the
inclusion morphism. The \emph{Kustin--Miller unprojection ring} $S$
of the pair $J \subset R$ is the quotient ring
\[
          S =  \frac {R[T]} {(Tu-\phi(u)  \bigm| u \in J) },
\]
where $T$ is a new variable. The ring $S$ is, up to isomorphism, 
independent of the choice of $\phi$.
The original definition of Kustin and Miller  \cite{KM}  was
using projective resolutions, compare
Subsection~\ref {subs!generalities_about_KM_complexes}  below.

\subsection {The Kustin--Miller complex construction }  \label{subs!generalities_about_KM_complexes}

The following construction, which is due to Kustin and Miller  \cite {KM},
 will be important in Section~\ref{sec!application_to_resolutions},
where we identify  the minimal graded free resolution of $k[\Delta (d,m)]$.

Assume $R$ is a polynomial ring over a field  with the degrees of all variables positive, 
 and $I \subset J \subset R$ are two 
homogeneous ideals of $R$ such
that both  quotient rings $R/I$ and $R/J$ are Gorenstein and $\dim R/J = \dim R/I -1$. 
We define $k_1,k_2 \in \Z$ such that $\omega_{R/I} = R/I(k_1)$   and  
$\omega_{R/J} = R/J(k_2)$, compare~\cite[Proposition~3.6.11]{BH}, and assume that $k_1 > k_2$. 
We fix a graded homomorphism $\phi \in \Hom_{R/I}(J,R/I)$ of degree $k_1-k_2$ 
such that $\Hom_{R/I}(J,R/I)$ is generated as an $R/I$-module by the subset 
$\{i, \phi \}$, where $i$ denotes the inclusion morphism,
compare~Subsection~\ref{subs!defnofkustinmillerunprojection}. We denote by 
$S = R[T]/Q$  the Kustin--Miller unprojection  ring of the pair 
$J \subset R/I$ defined by $\phi$, where $T$ is a new variable of degree $k_1-k_2$.
We have that $Q = (I, Tu-\phi(u)  \bigm| u \in J)$ and that $S$ is a graded algebra.

We denote by  $g= \dim R - \dim R/J$ the codimension of the ideal $J$ of $R$.
Let  
\[
  C_J : \quad    0 \to  R=A_{g} \to A_{g-1}  \to  \dots  \to A_1 \to R=A_0 
\] 
and
\[ 
  C_I : \quad  0 \to R= B_{g-1} \to  \dots  \to B_1 \to R = B_0 
\]
be the minimal graded free resolutions of $R/J$ and $R/I$ respectively as 
$R$-modules. Due to the Gorensteiness of $R/J$ and $R/I$ they are both
self-dual. 
We denote by  
$a_i \colon A_i \to A_{i-1}$ and 
$b_j \colon B_j \to B_{j-1}$ the differential maps. 
In the following, for an $R$-module $M$ we denoted 
by $M'$ the $R[T]$-module $M \otimes_R R[T]$. 

Kustin and Miller constructed in \cite {KM}
a graded free resolution $C_S$ of 
$S$ as $R[T]$-module of the form 
\[
   C_S : \quad   0 \to F_{g} \to F_{g-1} \to \dots  \to F_1 \to F_0 \to S\to 0,
\] 
where,  when $g \geq 3$,   
\begin {eqnarray*}
  &  &  F_0 =  B_0', \quad \quad  F_1 = B_1' \oplus A_1'(k_2-k_1),  \\\ 
  &  &   F_i = B_i' \oplus A_i'(k_2-k_1) \oplus  B_{i-1}'(k_2-k_1), \quad \quad   \text{for} \; \; 2 \leq i \leq g-2, \\
  &  &  F_{g-1} = A_{g-1}'(k_2-k_1) \oplus  B_{g-2}'(k_2-k_1), \quad \quad F_{g} =  B_{g-1}'(k_2-k_1), 
\end {eqnarray*}
cf. \cite [p.~307, Equation (3)]{KM}.  When $g=2$ we have 
\[ 
     F_0 =  B_0', \quad   F_1 =  A_1'(k_2-k_1), \quad  F_2 =  B_1'(k_2-k_1).
\]

We will now describe the differentials of the complex  $C_S$.
We denote the rank of the free $R$-module $A_1$ by $t_1$, since $C_J$ is
self-dual  $t_1$ is also the rank of the free $R$-module $A_{g-1}$.
We fix $R$-module bases $e_1, \dots , e_{t_1}$ of $A_1$  and 
$\hat{e}_1, \dots , \hat{e}_{t_1}$ of $A_{g-1}$. We define, for
$1 \leq i \leq t_1$,  $c_i,  \hat{c}_i  \in R$ by $ a_1 ( e_i) = c_i 1_R$  and 
$a_g(1_R) = \sum_{i=1}^{t_1} \hat{c}_i \hat{e}_i$. By Gorensteiness
we have that  $c_i, \hat{c}_i \in J$ for all $1 \leq i \leq t_1$.
For $1 \leq i \leq t_1$, let $l_i \in R$  be a lift in $R$ of $\phi (c_i) $ and
let $\hat{l}_i \in R$ be a lift in $R$ of $\phi (\hat{c}_i) $.
For an R-module $A$ we set $A^* =  \Hom_R (A,R)$. For an $R$-basis 
$f_1, \dots f_t$ of $A$ we denote by 
$f_1^*, \dots ,  f_t^*$ the basis of $A^*$ dual to it.

Denote by $\tilde{\alpha}^{d}_{g-1} \colon A_{g-1}^* \to R=B_{g-1}^*$ the $R$-homomorphism
with $\tilde{\alpha}^{d}_{g-1} (\hat{e}_{i}^*) = \hat{l}_i 1_R$ for $1 \leq i \leq t_1$. 
Taking into account the self-duality of $C_I,C_J$, 
we have that $\tilde{\alpha}^{d}_{g-1}$ extends to a 
chain map $\tilde{\alpha}^{d} \colon  C_J^* \to C_I^*$.  We denote by
$\tilde{\alpha} \colon C_I \to C_J$ the chain map dual to $\tilde{\alpha}^{d}$. The map 
$\tilde{\alpha}_0 \colon  B_0 = R \to R=A_0$ is multiplication by an invertible element, say $w$, 
of $R$, cf. \cite{P1},  and we set $\alpha = \tilde{\alpha} /w$.

We will now define a chain map $\beta \colon C_J \to C_I[-1]$. We first define
$ \beta_1 \colon A_1 \to R = B_0$ by $ \beta_1 ( e_i ) = -l_i 1_R$. 
We obtain a chain map $\beta \colon C_J \to C_I[-1]$ by extending 
$\beta_1$. Moreover,
$\beta_g \colon A_{g}=R \to R = B_{g-1}$ is multiplication by a nonzero
constant $u \in R$.   By \cite[p.~308]{KM} there exists a homotopy map  $h \colon C_I \to C_I$ 
with $h_0 \colon B_0 \to B_0$ and $h_{g-1} \colon B_{g-1} \to B_{g-1}$ being the zero 
maps and 
\[
       \beta_i \alpha_i = h_{i-1} b_i + b_i h_i,
\]
for $1 \leq i \leq g$.

Finally, following \cite[p.~307]{KM}, we have that the differential maps
$f_i \colon F_i \to F_{i-1}$ of the complex $C_S$ are given in block format by the
following formulas
\[
     f_1  = \left [   \begin{array}{cc}
                 b_1  &   \beta_1+Ta_1 
             \end{array} \right ],   \quad  \quad  \quad 
     f_2 =  \left [   \begin{array}{ccc}
                 b_2  &   \beta_2  &  h_1 + T I_1 \\
                 0    &    -a_2   &  -\alpha_1   
             \end{array} \right ],
\]

\[
    f_i =  \left [   \begin{array}{ccc}
                 b_i  &   \beta_i  &  h_{i-1} + (-1)^i T I_{i-1}  \\
                 0    &    -a_i   &  -\alpha_{i-1}  \\
                 0   &      0     &   b_{i-1}
             \end{array} \right ]    \quad  \quad  \text {for } 3 \leq i \leq g-2,
\]

\[  
    f_{g-1} = \left [   \begin{array}{cc}
                 \beta_{g-1}  &  h_{g-2} + (-1)^{g-1} T  I_{g-2} \\
                  -a_{g-1}    &   -\alpha_{g-2}  \\
                       0      &     b_{g-2}
             \end{array} \right ],
\]

\[  
    f_g = \left [   \begin{array}{c}
                 -\alpha_{g-1}  + (-1)^g u^{-1}Ta_g   \\
                    b_{g-1}  
             \end{array} \right ],
\]
where $I_t$ denotes the identity $\rank B_t \times \rank B_t$ matrix.

The resolution $C_S$ is, in general,
not minimal  \cite[Example~5.2]{BP}.  However, in the cases of  stacked 
and cyclic polytopes it is  minimal, see \cite{BP} and Theorem~\ref{thm!applicationtocyclicpolytopes}.
In the following we will call $C_S$  the \emph{Kustin--Miller complex construction}. 
We refer the reader to Subsection~\ref{subs!examplesforKMconstruction} for explicit
examples of this construction.

\section  { The main theorem for  $d$ even}  \label {sec!disscussion_for_d_even}
   
We fix a field $k$,  and assume that $d,m$ are integers with  $d$ even  and  $2 \leq d < m -1$.
(The case $m = d+1$ is discussed in Subsection~\ref {sub!casedevenandmisdplus}.)
We  set $a = (d+2)/2$, and denote by $k[\Delta (d,m)]$ the Stanley-Reisner ring
of the simplicial complex $\Delta (d,m)$.

The following lemma is an almost immediate consequence of Theorem~\ref {thm!facesofcyclicpolytopinfo}.

\begin {lemma}  \label {lemma!mainthmstanleyreisneridela}
    We have 
      \[
             k [ \Delta (d,m)]  \iso  k[x_1, \dots , x_m]/  (I_{a,1,m-1},  I_{a,2,m}).
      \]
\end {lemma}

\begin {proof} 

Denote by ${\mathcal A}$ the set of minimal monomial generators of the ideal $(I_{a,1,m-1}, I_{a,2,m})$. 
We first show that if $x_V \in {\mathcal A}$, then $V$ is not a face of $\Delta (d,m)$.
Assume  $x_V$ is a monomial generator of $I_{a,1,m-1}$, the case  $x_V$ is
a monomial generator of $I_{a,2,m}$ follows by the same arguments. Since $\# V = a$, 
we have that the number of odd
contiguous subsets of $V$ is at least $a-1$. Since $a-1 = d/2 > d/2 -1 = d -a$, by
Theorem~\ref {thm!facesofcyclicpolytopinfo} $V$ is not a face of $\Delta (d,m)$. 

Assume now  $W \subset \{ 1, \dots ,m \}$ is a subset with $\#W \leq d$.
We will show that if $W$ is  not a face of $\Delta (d,m)$
then there exists  a monomial generator $x_V \in {\mathcal A}$ with
$V \subset W$. By Theorem~\ref {thm!facesofcyclicpolytopinfo} $\# W \geq a$.
We will argue by induction on the cardinality of $W$.

Denote by $p$ the number of the odd contiguous subsets of $W$ considered as a 
subset of $\{1, \dots ,m\}$,
and, for $w \in W$, by $p_w$ the number of the odd contiguous subsets of $W \setminus \{w\}$
also considered as a subset of $\{1, \dots ,m\}$.
By  Theorem~\ref {thm!facesofcyclicpolytopinfo} $p > d - \#W$.
If $\#W = a$, then $p > d - \#W$ implies that $W$ has at least $d-a+1 = a-1 = \#W-1$ odd contiguous 
subsets,  and we set $V=W$.

Assume for the rest of the proof that $\# W > a$. 
By the inductive hypothesis it is enough to show
that there exists $w \in W$ such that $W \setminus \{w \}$ is not a face of $\Delta (d,m)$.
Hence, by  Theorem~\ref {thm!facesofcyclicpolytopinfo} it is enough to show that there 
exists $w \in W$  with $p_w > d- \#W + 1 $.

We call a nonempty  $X \subset W$ a gc-subset if there exist $i \leq j$ with $i-1 \notin W$, 
$j+1 \notin W$ such that $X=\{i,i+1, \dots , j \}$. It is obvious that a contiguous subset
of $W$ is a gc-subset, and that a gc-subset of $W$ is contiguous if and only if
contains neither $1$ nor $m$.

If $W$ contains
a gc-subset of even cardinality, say $\{ i, i+1, \dots , j \}$ we  set $w = m$ if $j = m$, while
if $j \not= m$ we set $w=i$. In the first case, since $i = 1$ contradicts $\#W \leq d$, we have 
that $p_w = p + 1$, so  $p_w > d- \#W +1 $ follows. Similarly,  for the second case again 
$p_w = p + 1$ and  $p_w > d- \#W +1 $ follows.

Assume for the rest of proof that all gc-subsets of $W$ are of odd cardinality.  First
assume that  $W$ contains a gc-subset $\{ i, i+1,  \dots , j \}$ 
of odd cardinality at least $3$, and  set $w = i+1$.
Since $(i,j) = (1,m)$ is impossible by $\#W \leq d$, 
it is clear that $p_w = p+1$, so again $p_w > d- \#W +1 $.

So we can assume for the rest of the proof that all  gc-subsets of $W$ are of 
cardinality  $1$. We either set $w = m$ if $m \in W$, or if $m \notin W$  we set 
$w$ to be the smallest element of $W$. 
If $m \in W$ and $1 \in W$ we have $p_w = p = \#W-2$, and $p > d - \#W$ implies 
$2\#W  - 2 > d$, so since $d$ is even $2\#W > d+3$, hence $p_w > d - \#W +1$.
If $m \in W$ and $1 \notin W$, we have $p_w = p = \#W-1$, and    $p_w > d- \#W +1 $
is equivalent to $2\#W > d+2$, which is true by the assumption  $ \#W > a = (d+2)/2$.
If $m \notin W$ and  $1 \in W$ the argument is exactly symmetric to the
case $m \in W$ and $1 \notin W$. If  $m \notin W$ and  $1 \notin W$, we have 
$p_w = p -1 = \#W -1 $ and  $p_w > d- \#W +1$ is equivalent to  $2\#W > d+2$, 
which is true by the assumption  $ \#W > a = (d+2)/2$.
This finishes the proof of 
Lemma~\ref {lemma!mainthmstanleyreisneridela}. \end {proof}

We now further assume that $d$ is an even integer with $d \geq 4$, the case $d=2$ is discussed in 
Subsection~\ref {sub!thecasedis2}.
We set $R = k[x_1, \dots ,x_m,z]$, where we put degree $1$ for all variables. 
We consider the ideals  $I = (I_{a,1,m-1},  I_{a,2,m})$   and $J=(I_{a-1, 2,m-1}, z  I_{a-2, 3, m-2})$
of $R$. (When we need to be more precise  we will
also use the notations $I_{d,m}$ for $I$ and $J_{d,m}$ for $J$.)
It is clear that $I \subset   (I_{a-1, 2,m-1})$, hence $I \subset J$. Moreover,
using Lemma~\ref{lemma!mainthmstanleyreisneridela},
$R/I \iso k[\Delta (d,m)][z]$ and $R/J \iso k[\Delta (d-2,m-1)][x_1,x_m]$. Consequently,
both rings $R/I$ and $R/J$ are Gorenstein by \cite [Corollary 5.6.5]{BH}, and $\dim R/J = \dim R/I -1$.

The  proof of the following key lemma 
will be given in Subsection~\ref{subsec!identificationofHommodule}.

\begin {lemma}  \label{lem!identificationofHpommodule} 
  There exists unique $\phi \in \Hom_{R/I}(J,R/I)$ such that
  $\phi (v) = 0$ for all $v \in I_{a-1,2,m-1}$ and $\phi (zw) = wx_1x_m$ for all
  $w \in I_{a-2,3,m-2}$.  Moreover, the $R/I$-module $\Hom_{R/I}(J,R/I)$ is
  generated by the set $\{ i, \phi \}$, where $i \colon J \to R/I$ denotes the
  inclusion homomorphism.
\end {lemma}

Taking into account Lemma~\ref{lem!identificationofHpommodule},
the Kustin--Miller unprojection ring $S$ of the pair $J \subset R/I$  is equal to 
\[
    S = \frac {(R/I)[T]}{(Tu-\phi (u) \bigm| u \in J) }.
\]
We  extend the grading of $R$ to a grading of $S$ by putting the degree of the new variable 
$T$ equal to $1$. By Lemma~\ref{lem!identificationofHpommodule}  $S$ is a graded $k$-algebra.
Our main result for the case $d$  even is the following theorem.

\begin {theorem}  \label {theorem!mainthm_cyclicpolytopes_deven}  
         The element $z \in S$ is $S$-regular, and there is 
         an isomorphism of graded $k$-algebras   
         \[
                        S/(z)  \iso   k[\Delta (d,m+1)].
         \]
\end {theorem}

\begin {proof} Denote by $Q \subset R[T]$ the ideal
\[    
      Q = (I,z) + (Tu-\phi (u) \bigm| u \in J)     \subset R[T].
\]
By the definition of $S$ we have $S/(z) \iso R[T]/Q$. By the definition of $\phi$
we have  $Q = (I_{a,1,m}, TI_{a-1,2,m-1},z)$.
Hence, Lemma~\ref {lemma!mainthmstanleyreisneridela}  implies that $S/(z) \iso   k[\Delta (d,m+1)]$.
As a consequence, $\dim S/(z) = \dim S -1$, and since by
\cite [Theorem~1.5]{PR}  $S$ is Gorenstein, hence Cohen--Macaulay, we get that $z$ is
 $S$-regular.  \end {proof}

\begin {example}  Assume $d=4$ and $m = 6$. We have 
\[
    I=  (x_2x_4x_6,x_1x_3x_5)  ,  \quad \quad  J = (x_2x_4, x_2x_5,  x_3x_5,   zx_3,   zx_4 )
\]
and  
{ \small
\[
    S = k[x_1, \dots ,x_6,T,z]/(I , Tx_2x_4, Tx_2x_5,  Tx_3x_5, x_3 (zT- x_1x_6) ,        x_4 (zT- x_1x_6)).
\]
}
\end {example}

\subsection  { The case $d=2$ and $d+1 < m$ }  \label {sub!thecasedis2}

Assume $d=2$ and $d+1 < m$. It is clear that $\Delta (d,m)$ is just the (unique) triangulation of the
$1$-sphere $S^1$ having $m$ vertices. Hence $\Delta (d,m+1)$ is a stellar subdivision
of $\Delta (d,m)$, and the results of \cite {BP} apply.

In more detail, set $R=  k[x_1, \dots ,x_m,z]$, with the degree of all
variables equal to $1$. Consider the ideals  $I = (I_{2,1,m-1} , I_{2,2,m})$ and
$J = (I_{1,2,m-1}, z)$ of $R$.  (When we need to be more precise we will
also use the notations $I_{2,m}$ for $I$ and $J_{2,m}$ for $J$.)
Clearly $k[\Delta (d,m)][z] \iso  R/I$.
Moreover, we have that $I \subset J$, that $J \subset R/I$ is a codimension $1$ ideal of $R/I$
with $R/J$ Gorenstein, and that if we denote by $S$ the Kustin--Miller unprojection ring
of the pair $J \subset R/I$ we have $S/(z) \iso k[\Delta (d,m+1)]$.
Moreover, arguing as in the proof of Theorem~\ref{theorem!mainthm_cyclicpolytopes_deven}   
we get that $z$ is an $S$-regular element.

\subsection  { The case $d$ is even and  $m=d+1$ }  \label {sub!casedevenandmisdplus}

Assume $d \geq 2$ is even and $m=d+1$.  We have that
\[
   k[\Delta (d,m)] \iso k[x_1, \dots , x_{m}]/( \prod_{i=1}^{d+1} x_i)
\]
and
\[
     k[\Delta (d,m+1)] \iso k[x_1, \dots , x_{m+1}]/(\prod_{i=0}^{d/2} x_{2i+1}, \prod_{i=1}^{(d/2)+1} x_{2i}).
\]

We set $R= k[x_1, \dots ,x_m,z]$, with the degree of all variables equal to $1$.
Consider the ideals $I=(\prod_{i=1}^{d+1} x_i)$ and 
$J = (\prod_{i=1}^{d/2} x_{2i}, z\prod_{i=1}^{(d/2)-1} x_{2i+1})$
of $R$. (When we need to be more precise we will
also use the notations $I_{d,m}$ for $I$ and $J_{d,m}$ for $J$.)
We have $I \subset J$, that $J \subset R/I$ is a codimension $1$ ideal of $R/I$
with $R/J$ Gorenstein, and that if we denote by $S$ the Kustin--Miller unprojection
ring of the pair $J \subset R/I$ we have
$ S/(z) \iso k[\Delta (d,m+1)]$.
Moreover, arguing as in the proof of Theorem~\ref{theorem!mainthm_cyclicpolytopes_deven}   
we get that $z$ is an $S$-regular element.

\subsection {Proof of Lemma~\ref{lem!identificationofHpommodule} } \label{subsec!identificationofHommodule}

We start the proof of Lemma~\ref{lem!identificationofHpommodule}. Recall that 
$I = (I_{a,1,m-1},  I_{a,2,m})$   and $J=(I_{a-1, 2,m-1}, z  I_{a-2, 3, m-2})$.
Since $J$ is a codimension $1$ ideal of $R/I$ and $R/I$ is Gorenstein, hence Cohen--Macaulay,
there exists $b \in J$ which is $R/I$-regular. Write $b = b_1 + z b_2$, with $b_1 \in I_{a-1,2,m-1}^e$ and
 $b_2 \in I_{a-2,3,m-2}^e$, where $I_*^e$ denotes the ideal of $R/I$ generated by $I_*$.
Consider the element 
\[
    s_0 = \frac {b_2 x_1x_m} {b} \in K(R/I), 
\]
where $K(R/I)$ denotes the total quotient ring of $R/I$, 
that is the localization
of $R/I$ with respect to the multiplicatively closed subset of regular elements of $R/I$,
cf.~\cite[p.~60]{Ei}.   We need the following lemma.

\begin {lemma}  \label {lemma!basicidentitiesforphi24}
   (a)  We have that  $x_1x_mvw = 0$ (equality in $R/I$) for all $v \in I_{a-1,2,m-1}$ and 
        $w \in I_{a-2,3,m-2}$.  

  (b)  We have $s_0 zw = w x_1x_m$ (equality in $K(R/I)$) for all $w \in I_{a-2,3,m-2}$.  
\end {lemma}

\begin {proof}  
  Proof of (a).  It  is enough to show  that $x_1x_mx_Vx_W = 0$ in $k[\Delta (d,m)]$, 
  whenever $x_V$ is a generating monomial of $I_{a-1,2,m-1}$
  and  $x_W$ is a generating monomial of $I_{a-2,3,m-2}$,
  with  $V \subset \{2, \dots ,m-1 \}$ and  $W \subset \{3, \dots ,m-2 \}$.
    Consider the set
  $A = \{1,m\} \cup V \cup W$.  If $2 \notin V$ it is clear that $x_1x_V = 0$  and, similarly,
  if $m-1 \notin V$ we have $x_mx_V = 0$. 

  Hence for the rest of the proof we can assume that $2 \in V$ and $m-1 \in V$.  Denote by $A_1=\{1, \dots , p\}$ the initial 
 segment of $A$, and by $A_2$ the final segment of $A$.  Since $2,m-1 \notin W$, we necessarily have
 that all odd elements of $A_1 \setminus \{ 1 \}$  are in $W \setminus V$, and all even elements of
 $A_1$ are in $V \setminus W$. If the largest element $p$ of $A_1$
 is not in $V$,   the monomial with support   $(V \setminus A_1) \cup \{1,3, \dots,  p \}$ is in $I$, hence 
 $x_1x_Vx_W = 0$. By a similar argument, if the smallest element  of $A_2$ is not in $V$ we get $x_mx_Vx_W = 0$. 
 So we can assume that both the largest element of   $A_1$  and the smallest element of 
 $A_2$ are in $V$.  By the above discussion, this implies that $ \#(A_1 \cap V) = \#(A_1 \cap W)+1 $ and 
 $\#(A_2 \cap V) = \#(A_2 \cap W)+1$, hence  $\# W_a  =  \# V_a +1$,
 where we set $V_a = V \setminus (A_1 \cup A_2)$ and $W_a = W \setminus (A_1 \cup A_2)$.
 Hence there exists a 
 contiguous subset of $V_a \cup W_a$, say  $A_3 = \{i, i+1, \dots ,j \}$,
 which starts with an element of $W \setminus V$ then either stops or continuous
 with an element of $V \setminus W$ and finally finishes with an element of $W \setminus V$. The monomial
 with support in  $(V \setminus A_3) \cup \{i,i+2, \dots,  j\}$  is in $I$, hence we get
 $x_Vx_W = 0$ which finishes the proof of part (a) of Lemma~\ref{lemma!basicidentitiesforphi24}. 

  We now prove part (b) of the lemma.  It is enough to show that 
  $(b_1+zb_2)wx_1x_m = zw(b_2 x_1 x_m)$, for all $ w \in W$. For
  that it is enough to show $x_1x_mb_1w = 0$, which follows from part (a).
\end {proof}

Using Lemma~\ref{lemma!basicidentitiesforphi24}, multiplication by $s_0$, which a priori is only
an $R/I$-homomorphism $R/I \to K(R/I)$, maps $J$ inside $R/I$, so defines an 
$R/I$-homomorphism $\phi \colon J \to R/I$. By the same Lemma~\ref{lemma!basicidentitiesforphi24},
 we have that  $\phi (v) = 0$, for all $v \in I_{a-1,2,m-1}$, and $\phi (zw) = wx_1x_m$, for all
  $w \in I_{a-2,3,m-2}$.  Since an $R/I$-homomorphism is uniquely determined
by its values on a generating set, the uniqueness of $\phi$ stated in 
Lemma~\ref{lem!identificationofHpommodule} follows.

We will now prove the part of  Lemma~\ref{lem!identificationofHpommodule}  stating that the 
$R/I$-module $\Hom_{R/I}(J,R/I)$ is generated by the set $\{ i, \phi \}$.  By the arguments 
contained in the proof of \cite[Theorem~5.6.2] {BH}, we have isomorphisms
\[
  \omega_{k[\Delta (d,m)]} \iso k[\Delta (d,m)](0), \
       \quad  \omega_{k[\Delta (d-2,m-1)]} \iso k[\Delta (d-2,m-1)](0),
\]
of graded $k$-algebras, where $\omega_R$ denotes the canonical $R$-module.
Consequently,  since $R/I \iso k[\Delta (d,m)][z]$, 
$R/J \iso k[\Delta (d-2,m-1)][x_1,x_m]$ we get
\begin{equation}  \label {eqn!aboutomegasofroveriroverj}
   \omega_{R/I}  \iso (R/I)(-1) \quad \text { and }  \quad \omega_{R/J}  \iso (R/J)(-2).
\end {equation}
Combining  (\ref {eqn!aboutomegasofroveriroverj})  with the short exact sequence 
(\cite[p.~563] {PR})
\[
   0 \to \omega_{R/I} \to \Hom_{R/I}  (J,\omega_{R/I}) \to \omega_{R/J} \to 0,
\]
we get the short exact sequence
\[
   0 \to  R/I \to \Hom_{R/I}  (J,R/I) \to (R/J)(-1) \to 0.
\]
As a consequence, $\Hom_{R/I}  (J,R/I)$ is generated as an $R/I$-module by
the subset $\{ i, \psi \}$, whenever $\psi \in  \Hom_{R/I}  (J,R/I) $ has homogeneous degree
$1$ and  is not contained in the $R/I$-submodule of  $\Hom_{R/I}  (J,R/I)$ generated
by the inclusion homomorphism $i$.  Hence, to prove $\Hom_{R/I}  (J,R/I) = (i, \phi)$
is enough to show that there is no $c \in R/I$ with $\phi = c i$.
Assume such $c$ exists. Let $w \in I_{a-2,3,m-2}$ be a fixed monomial generator. 
We then have $czw =  \phi(zw) = wx_1x_m$ (equality in $R/I$), and since $R/I$
is a polynomial ring with respect to $z$ we get $wx_1x_m =0$, which is impossible,
since $I=(I_{a,1,m-1}, I_{a,2,m})$.  Hence
 $\Hom_{R/I}  (J,R/I) = (i, \phi)$, which finishes the
proof of Lemma~\ref{lem!identificationofHpommodule}.

\section  { The main theorem for  $d$ odd}  \label {sec!disscussion_for_d_odd}

Assume $k$ is a fixed field, and $d,m$ two integers with $d$ odd and $5 \leq d < m-1$,  the cases $d=3$
and $m = d+1$ are discussed in Subsections~\ref{sub!casedis3and} and \ref{sub!casedoddandmisdplusone}
respectively. We set $a = (d+1)/2$. Combining  
Proposition~\ref {lemma!mainthmstanleyreisneridela} 
with \cite[Exerc.~5.2.18]{BH} we get the following proposition.

\begin {prop} \label{prop!reduction_step_from_odd_d_to_d_1}  We have
\[
     k[\Delta (d,m)]  \ \iso  k[x_1, \dots ,x_m]/ (I_{a,2,m-1}, x_1x_m I_{a-1,3,m-2}).
\]
\end {prop}

\begin {remark}  \label{rem!remark_about_reduction_from_odd_d_to_d_1}
By Proposition~\ref{prop!reduction_step_from_odd_d_to_d_1}   and \cite[Exerc.~5.2.18]{BH},  for $d \geq 5$ odd
the ideal defining   $k[\Delta (d,m)]$ is related to the ideal 
defining  $k[\Delta (d-1,m-1)]$.  We will use this in what follows
to reduce questions for  $d$ odd to the easier case $d$ even.
A similar remark also applies when $d =3$.
\end {remark}

We set $R = k[x_1, \dots ,x_m,z_1,z_2]$, where we put degree $1$ for all variables.
Consider the ideals  $I=(I_{a,2,m-1}, x_1x_m I_{a-1,3,m-2})$ and 
$J= (I_{a-1, 2,m-2}, $\\ 
$z_1z_2  I_{a-2, 3, m-3})$
of $R$.
It is clear that $I \subset  (I_{a-1, 2,m-2})$, hence $I \subset J$. By
Proposition~\ref{prop!reduction_step_from_odd_d_to_d_1} we have that
$R/I \iso k[\Delta (d,m)][z_1,z_2]$ and 
$R/J \iso k[\Delta (d-2,m-1)][x_1,x_{m-1}, x_{m}]$. Consequently,
both rings $R/I$ and $R/J$ are Gorenstein by \cite [Corollary 5.6.5]{BH},  
and $\dim R/J = \dim R/I -1$. The following lemma is the analogue of 
Lemma~\ref{lem!identificationofHpommodule} for the case $d$  odd.

\begin {lemma}  \label{lem!identificationofHpommodulecasedodd} 
  There exists unique $\phi \in \Hom_{R/I}(J,R/I)$ such that
  $\phi (v) = 0$ for all $v \in I_{a-1,2,m-2}$ and $\phi (z_1z_2w) = x_1x_{m-1}x_{m}w$ for all
  $w \in I_{a-2,3,m-3}$.  Moreover, the $R/I$-module $\Hom_{R/I}(J,R/I)$ is
  generated by the set $\{ i, \phi \}$, where $i \colon J \to R/I$ denotes the
  inclusion homomorphism.
\end {lemma}

\begin {proof}
  Taking into account Proposition~\ref{prop!reduction_step_from_odd_d_to_d_1}
   and Remark~\ref{rem!remark_about_reduction_from_odd_d_to_d_1},
  Lemma~\ref{lem!identificationofHpommodulecasedodd} 
  follows  by the same arguments as  Lemma~\ref{lem!identificationofHpommodule}.
\end {proof}

Taking into account Lemma~\ref{lem!identificationofHpommodulecasedodd},
the Kustin--Miller unprojection ring $S$ of the pair $J \subset R/I$ is equal to
\[
      S = \frac {(R/I)[T]}{(Tu-\phi (u) \bigm| u \in J) }.
\]
We  extend the grading of $R$ to a grading of $S$ by putting the degree of the new variable 
$T$ equal to $1$.
Lemma~\ref{lem!identificationofHpommodulecasedodd} tells us  that $S$ is a graded $k$-algebra.
Our main result for the case $d$ odd is the following theorem.

\begin {theorem}  \label {theorem!mainthm_cyclicpolytopes_dodd}  
         The sequence $z_1,z_2 \in S$ is $S$-regular, and there is 
         an isomorphism of graded $k$-algebras   
         \[
                        S/(z_1,z_2)  \iso   k[\Delta (d,m+1)].
         \]
\end {theorem}

\begin {proof}  
Denote by $Q \subset R[T]$ the ideal
\[
        Q = (I,z_1,z_2) + (Tu-\phi (u) \bigm| u \in J)     \subset R[T].
\] 
By the definition of $S$ we have $S/(z_1,z_2) \iso R[T]/Q$.

Denote by $g \colon R[T] \to  R[x_{m+1}]$ the  $k$-algebra 
isomorphism  which is uniquely specified by 
$g(z_i) = z_i$ for $i=1,2$, $g(x_i) = x_i$ for $1 \leq i \leq m-1$, $g(x_m) = x_{m+1}$ and 
$g(T) = x_m$. It is easy to see that  $g (Q) = (I_{d,m+1},z_1,z_2)$.
Since $g$  is an  isomorphism, we have using Proposition~\ref{prop!reduction_step_from_odd_d_to_d_1} 
that 
\[
     R[T]/Q \iso   R[x_{m+1}]/ (I_{d,m+1},z_1,z_2) \iso  k[\Delta (d,m+1)],
\]
hence 
$S/(z_1,z_2)  \iso   k[\Delta (d,m+1)]$.  As a consequence, $\dim S/(z_1,z_2) = \dim S -2$, 
and since by \cite [Theorem~1.5]{PR}  $S$ is Gorenstein, hence Cohen--Macaulay, we get that 
$z_1,z_2$ is an $S$-regular sequence.  \end {proof}

\subsection  { The case $d=3 $ and $d +1 < m$ }  \label{sub!casedis3and}

Assume $d=3$ and $d +1 < m$. Combining \cite[p.~229, Exerc.~5.2.18]{BH} 
with the discussion of Subsection~\ref {sub!thecasedis2} we have the following picture.
Set $R=  k[x_1, \dots ,x_m,z_1,z_2]$,   where we put degree $1$ for all variables.
Consider the ideals  $I= (I_{2,2,m-1}, x_1x_m I_{1,3,m-2})$ and 
 $J= (I_{1,2,m-2},z_1z_2)$ of $R$. 
Then  $ k[\Delta (d,m)][z_1,z_2] \iso  R/I $.
Moreover, we have $I \subset J$, that $J \subset R/I$ is a codimension $1$ ideal of $R/I$
with $R/J$ Gorenstein, and that if we denote by $S$ the Kustin--Miller unprojection
ring  of the pair $J \subset R/I$
then $z_1,z_2$ is an $S$-regular sequence and $ S/(z_1,z_2) \iso k[\Delta (d,m+1)].$

\subsection  { The case $d$ is odd and  $m=d+1$ } \label{sub!casedoddandmisdplusone}

Assume $d \geq 3$ is odd and $m=d+1$.  We  have 
\[
   k[\Delta (d,m)] \iso k[x_1, \dots , x_{m}]/(\prod_{i=1}^{d+1} x_i)
\]
and
\[
     k[\Delta (d,m+1)] \iso k[x_1, \dots , x_{m+1}]/(\prod_{i=0}^{(d+1)/2} x_{2i+1}, \prod_{i=1}^{(d+1)/2} x_{2i}).
\]
Set $R= k[x_1, \dots ,x_m,z_1,z_2]$, where we put degree $1$ for all variables.
Consider the ideals  $I = (\prod_{i=1}^{d+1} x_i)$ and 
$J= ( \prod_{i=1}^{(d+1)/2} x_{2i}, z_1z_2 \prod_{i=1}^{(d-1)/2} x_{2i+1})$
of $R$. 
We have $I \subset J$, that $J \subset R/I$ is a codimension $1$ ideal of $R/I$
with $R/J$ Gorenstein, and that if we denote by $S$ the Kustin--Miller unprojection
ring of the pair $J \subset R/I$ then $z_1,z_2$ is an $S$-regular sequence and 
$ S/(z_1,z_2) \iso k[\Delta (d,m+1)]$.

\section {Combinatorial interpretation of our construction} \label{sec!combinatorial_interpretation}

We fix  $d \geq 2$ even  and  $m \geq d+1$, and  we will give a combinatorial interpretation of the 
constructions of Section~\ref{sec!disscussion_for_d_even}.  We  introduce the notation
$R_{(m)} = k[x_1, \dots ,x_m,z]$. 
Consider the ideals $I_{d,m}$ and $J_{d,m}$ of $R_{(m)}$ 
as defined in  Section~\ref{sec!disscussion_for_d_even} if $d \geq 4$ and $m \geq d+2$, 
as defined in Subsection~\ref{sub!thecasedis2} if $d=2$ and $m \geq d+2$,  and as defined in Subsection~\ref{sub!casedevenandmisdplus}
if $d \geq 2$ and $m= d+1$.

Note that $I_{d,m}$ is the Stanley--Reisner ideal of $\Delta (d,m)$. We will inductively identify
$J_{d,m}$.   We set   
$P_{d,m} = I_{d,m} : (x_1x_m)$, then 
\[ 
   P_{d,m}  = I_{\romstar_{\Delta(d,m)}(\{1,m\})} + 
        (x_i \bigm| i \text{ is not a vertex of } \romstar_{\Delta(d,m)}(\{1,m\})),
\]
It is clear that the ideal $P_{d,m}$ of $R_{(m)}$ is monomial, and that no minimal
monomial generator of it involves the variables $x_1,x_m$ and $z$. 
We denote by $\hat{P}_{d,m}$ the ideal of $k[x_2, \dots x_{m-1},z]$ which has the
same minimal monomial generating set.

If $d= 2$ we have $J_{d,m} = (P_{d,m} , z)$. Assume
now $d \geq 4$.  It is easy to see that  the ideal $\hat{P}_{d,m}$
is equal to the image of the ideal  $I_{d-2,m-2}$ of $R_{(m-2)}$
under the $k$-algebra isomorphism $R_{(m-2)} \to k[x_2, \dots x_{m-1},z]$ 
that sends $z$ to $z$ and $x_i$ to $x_{i+1}$ for $1 \leq i \leq m-2$, 
hence  $\hat{P}_{d,m}$ is the Stanley--Reisner ideal of a simplicial complex isomorphic
to $\Delta (d-2,m-2)$.  The unprojection constructions described in 
Section~\ref{sec!disscussion_for_d_even}
and Subsections~\ref{sub!thecasedis2}, \ref{sub!casedevenandmisdplus} allow us to
pass from  the ideal $I_{d-2,m-2}$ of $R_{(m-2)}$ to the ideal
$I_{d-2,m-1}$ of $R_{(m-1)}$, which is the Stanley--Reisner ideal of 
$\Delta (d-2,m-1)$.   Denote by 
$Q_{d,m} \subset k[x_2, \dots ,x_m,z]$ the image of the ideal 
$I_{d-2,m-1}$ under the $k$-algebra isomorphism 
$R_{(m-1)} \to k[x_2, \dots ,x_m,z]$ that sends   $z$ to $x_m$,
$x_i$ to $x_{i+1}$ for $1 \leq i \leq m-2$,
and $x_{m-1}$ to $z$. It is then easy to see that  $J_{d,m}$ is
the ideal of  $R_{(m)}$ generated by the image 
of $Q_{d,m}$ under the inclusion of $k$-algebras
$k[x_2, \dots ,x_m,z] \to R_{(m)}$. In particular,
$R_{(m)}/(J_{d,m},x_1,x_m) \iso k[\Delta (d-2,m-1)]$,
as already observed above.

Assume now $d \geq 3$ is odd and $m \geq d+1$.   Consider the ideal 
$J$ as defined in  Section~\ref{sec!disscussion_for_d_odd}. 
Using Remark~\ref{rem!remark_about_reduction_from_odd_d_to_d_1},
a similar combinatorial interpretation  exists  for $J$ 
in terms of the $\Delta (d-2,m-2)$  related to the star of the face 
$\{1,m \}$  of  $\Delta (d,m)$ when $d \geq 5$, and  an
analogous statement when $d=3$. We leave the precise formulations to the reader.

\section {The minimal resolution of cyclic polytopes}  \label{sec!application_to_resolutions}

Combining the results of Sections~\ref{sec!disscussion_for_d_even}
and~\ref{sec!disscussion_for_d_odd}, we have that for
$d \geq 4$ and $d+1 < m$,  the Stanley-Reisner ring  $k[\Delta (d,m+1)]$
can be constructed from the Stanley--Reisner rings $k[\Delta (d,m)]$ and 
$k[\Delta (d-2,m-1)]$ using Kustin--Miller unprojection.
Moreover, we showed that
a similar statement is true also for the cases $d=2,3$ and $m=d+1$. 
Using the Kustin--Miller complex construction discussed 
in Subsection~\ref{subs!generalities_about_KM_complexes},
we can inductively build a graded free resolution of $S$,  hence
using  Proposition~\ref{prop!dividing_resolutions_by_graded_elements} below
of $k[\Delta (d,m+1)]$,  starting from the minimal graded 
free resolutions of $k[\Delta (d,m)]$ and $k[\Delta (d-2,m-1)]$.  The following
theorem, which will be  proven in Subsection~\ref{subs!pfofthemapplicationto}, tells us 
that in this way we get  a minimal resolution.   
Subsection~\ref{subs!examplesforKMconstruction} contains examples
demonstrating the theorem and a link to related 
computer algebra code.

\begin {theorem}  \label{thm!applicationtocyclicpolytopes}
For $d \geq 4$ and $d+1 < m$,  the graded free resolution of 
$k[\Delta (d,m+1)]$ obtained from the minimal graded free resolutions of 
$k[\Delta (d,m)]$ and $k[\Delta (d-2,m-1)]$
using the Kustin--Miller complex construction is minimal.
For $d =2$ or $3$  and $d+1 < m$, the graded free 
resolution of  $k[\Delta (d,m+1)]$ obtained  from the  minimal graded free resolution 
of $k[\Delta (d,m)]$ and the appropriate Koszul 
complex  (see~Subsections~\ref{sub!thecasedis2} and \ref{sub!casedis3and})
using the Kustin--Miller complex construction
is also minimal.
\end {theorem}

We remark that in the proof of Theorem~\ref{thm!applicationtocyclicpolytopes} we
do not use the calculation of the graded Betti numbers of $k[\Delta (d,m)]$ obtained by
Schenzel \cite{S}  for even $d$, and by Terai and Hibi \cite{TH} for odd $d$.
Not only that, but in Proposition~\ref{prop!keytechnicalforapplication}
we recover their results, without using Hochster's formula  or Alexander duality.

\subsection {Proof of Theorem~\ref{thm!applicationtocyclicpolytopes}} \label{subs!pfofthemapplicationto}

For the proof of Theorem~\ref{thm!applicationtocyclicpolytopes} we will need the following
combinatorial discussion. 

Assume $d \geq 3$ is odd,  $d+1 < m$  and  $1 \leq i \leq m-d-1$.
We set
\[
   \eta (d,m,i) = \binom  {m-[d/2]-2} {[d/2]+i} \binom { [d/2] + i-1} { [d/2]},
\]
compare \cite[p.~291]{TH}. We also set $\eta (d,m,0) = \eta(d,m,m-d) = 0$.

\begin {prop}   \label {prop!combinatorialformula}
 We have,  for $1 \leq i \leq  m-d$, 
\begin {equation}  \label {eqn!inductiveeqnforeta}
      \eta (d,m+1,i) =  \eta (d,m,i) + \eta (d,m,i-1) + \eta (d-2,m-1,i).
\end{equation}
(By our conventions, for $i=1$ the equality becomes  $\eta (d,m+1,1) =  \eta (d,m,1) + \eta (d-2,m-1,1)$,
while for $i = m-d$  it becomes  $\eta (d,m+1,m-d) =  \eta (d-2,m-1,m-d) + \eta (d,m,m-d-1)$.)
\end {prop}

\begin {proof}  Assume first $2 \leq i \leq  m-d-1$.  We will use twice the Pascal triangle identity
$\binom{k}{d} = \binom{k-1}{d} + \binom {k-1}{d-1}$. We have
{ \footnotesize

\begin {eqnarray*}
   \eta (d,m+1,i) &  = &  \binom  {m+1-[d/2]-2} {[d/2]+i} \binom { [d/2] + i-1} { [d/2]}   \\
  &  = & \left( \binom  {m-[d/2]-2} {[d/2]+i}+\binom  {m-[d/2]-2} {[d/2]+i-1} \right) \binom { [d/2] + i-1} { [d/2]}  \\
   & = &  \binom  {m-[d/2]-2} {[d/2]+i} \binom { [d/2] + i-1} { [d/2]}  + 
               \binom  {m-[d/2]-2} {[d/2]+i-1} \binom { [d/2] + i-1} { [d/2]}   \\
   & = & \ \eta (d,m,i) +  \binom  {m-[d/2]-2} {[d/2]+i-1} 
                  \left( \binom { [d/2] + i-2} { [d/2]} +  \binom { [d/2] + i-2} { [d/2]-1}   \right) \\
   & =  &  \eta (d,m,i) + \eta (d,m,i-1) + \eta (d-2,m-1,i).
\end{eqnarray*}
}
The special cases $i=1$ and  $i = m-d$ are proven by the same argument. \end {proof}

For the proof of Theorem~\ref{thm!applicationtocyclicpolytopes} we will also need the following 
general propositions, the first of which is well-known.

\begin {prop}  \label{prop!dividing_resolutions_by_graded_elements}
{\rm (\cite [Proposition~1.1.5] {BH}).}
Assume   $R= k[x_1, \dots ,x_n]$ is a polynomial
ring over a field $k$ with the degrees of all variables positive, and $I \subset R$ a homogeneous ideal.
Moreover, assume that $x_n$ is $R/I$-regular.  Denote by $cF$ the minimal graded free resolution
of $R/I$ as $R$-module. We then have that $cF \otimes_R R/(x_n)$ is the minimal graded
free resolution of   $R/(I,x_n)$  as  $k[x_1, \dots ,x_{n-1}]$-module, where we
used the natural isomorphisms   $ R \otimes_R R/(x_n) \iso R/(x_n) \iso k[x_1, \dots ,x_{n-1}]$.
\end {prop}

The following proposition is an immediate consequence of  Equation~(\ref{eqn!total_betti_as_tor_dimensions}).

\begin {prop}  \label {prop!sum_of_betti_remains_the_same}
Assume $k$ is a field and $ R_1 = k[x_1, \dots , x_n ]$, $R_2 = k[y_1, \dots ,y_n]$ are 
two polynomial rings with the degrees of all variables positive. 
Assume $I_1 \subset R_1$ is a monomial ideal, and
denote by   $I_2$  the ideal of $R_2$  generated by the image of  $I_1$ 
under the  $k$-algebra homomorphism $R_1 \to  R_2, \; x_i \mapsto y_i$, for $1 \leq i \leq n$.
Obviously $I_2$ is a homogeneous ideal of $R_2$.
We claim that for all $i \geq 0$ we have $b_i(R_2/I_2) = b_i (R_1/I_1)$
(of course the graded Betti numbers $b_{ij}$ of $R_2/I_2$ and $R_1/I_1$ 
may differ). \end {prop}

\begin {prop}  \label {prop!about_thesituation_replacing_t_by_t1_times_t2}
Assume $k$ is a field,   $R_1 = k[x_1, \dots   ,x_n,T ]$ and $R_2 = $ 
$k[y_1, \dots ,y_n,T_1,T_2]$
are two polynomial rings with the degrees of all 
variables positive, $\deg x_i = \deg y_i$, for $1 \leq i \leq n$,
and $\deg T = \deg T_1 + \deg T_2$. Assume $I_1 \subset R_1$ is a homogeneous ideal, and denote
by $I_2 \subset R_2$ the ideal generated by the image of $I_1$ under the graded $k$-algebra 
homomorphism $\phi \colon R_1 \to R_2$ specified by $\phi(x_i)= y_i$, for $1 \leq i \leq t$,   and 
$\phi(T) = T_1T_2$.
Denote by $cF_1$ the minimal graded free resolution of $R_1/I_1$ as $R_1$-module. Then $I_2$ is
a homogeneous ideal $R_2$, and the complex $cF_1 \otimes_{R_1} R_2$ is 
a minimal graded free resolution of $R_2/I_2$ as $R_2$-module. In particular, the 
corresponding graded Betti numbers $b_{ij}$ of $R_1/I_1$ and $R_2/I_2$ are equal.
\end {prop}

\begin {proof}   It is clear that $I_2$ is a homogeneous ideal of $R_2$.
 By \cite [Theorem~18.16]{Ei}  $\phi$ is flat. As a consequence,
 \cite [Proposition~6.1]{Ei} implies that the natural map $I_1 \otimes_{R_1} R_2 \to I_2$ is 
an isomorphism of graded $R_2$-modules.
By flatness, tensoring  the minimal graded free resolution of $I_1$ as $R_1$-module with  $R_2$
we get the minimal graded free resolution of $I_2$ as $R_2$-module, and 
Proposition~\ref{prop!about_thesituation_replacing_t_by_t1_times_t2} follows.
\end {proof}

Theorem~\ref{thm!applicationtocyclicpolytopes} will follow from the following 
more precise statement. Notice that, as we already mentioned before, the statements about the 
graded Betti  numbers
have been proven before by different arguments in \cite{S,TH}, but we do not
need to use their results.

\begin {prop}  \label {prop!keytechnicalforapplication}
Assume $d \geq 2$ and $d+1 < m$.  Set  $b_{ij} = b_{ij}(k[\Delta (d,m)])$.
Then the statement of Theorem~\ref{thm!applicationtocyclicpolytopes} is true for  $(d,m)$.
Moreover, we have that if $d$ is even then 
$b_{ij}= 1$ for $(i,j) \in \{ (0,0), (m-d,m) \}$, 
\[
   b_{i,d/2+i} = \eta (d+1,m+1, i) + \eta (d+1,m+1, m-d-i),
\]
for $1 \leq i \leq m-d-1$, and $b_{ij} = 0$ otherwise.  If $d$ is odd, then
$b_{ij}= 1$ for $(i,j)  \in \{ (0,0), (m-d,m) \}$, 
\[
   b_{i,[d/2]+i} = \eta (d,m, i), \quad b_{i,[d/2]+i+1} = \eta (d,m, m-d-i),
\]
for $1 \leq i \leq m-d-1$, and $b_{ij} = 0$ otherwise. 
 \end {prop}

\begin {proof}  We use induction on $d$ and $m$. 
If $d \geq 2$ and $m = d+2$ then $k[\Delta (d,m)]$ is a codimension $2$ 
complete intersection and everything is clear.  

The next step, is to notice that, for $d=2$ and $m \geq 3$,
Proposition~\ref{prop!keytechnicalforapplication} follows from 
\cite[Proposition~5.7]{BP}, since  $\Delta (2,m)$ is equal 
to $\Delta P_2(m)$ defined in \cite[Section~5]{BP}.

Now assume that $d$ is even with $d \geq 4$ and $d+3 \leq m$, and, by the 
inductive hypothesis,  Proposition~\ref {prop!keytechnicalforapplication}
holds for the values $(d-2,m-1)$ and $(d,m)$. An easy computation, 
taking into account Proposition~\ref {prop!combinatorialformula},
shows that the Kustin--Miller complex construction resolving $k[\Delta (d,m+1)]$
has the conjectured graded Betti numbers. Since no degree $0$ morphisms appear
it is  necessarily minimal.  This finishes the proof for $d$ even.

Assume now $d \geq 3$ is odd.  Combining \cite[Exerc.~5.2.18]{BH}  with 
Propositions \ref {prop!sum_of_betti_remains_the_same} and  \ref {prop!about_thesituation_replacing_t_by_t1_times_t2} 
 we get that, for $0 \leq i \leq m-d$,
\begin {equation} \label{eqn!keyequalityforbi_from_odd_d_to_d-1}
    b_i (k[\Delta (d,m)]) = b_i (k[\Delta (d-1,m-1)]).
\end{equation}
(Of course the graded Betti numbers $b_{ij}$ can, and in fact are, different for $k[\Delta (d,m)]$ and 
$k[\Delta (d-1,m-1)]$.)
So we can reduce the case $d$ odd to the case $d-1$, by doing an almost identical induction
on $(d,m)$ as in the case $(d-1,m-1)$, noticing that  the Kustin--Miller complex construction for 
$k[\Delta (d,m+1)]$  has to be minimal, since we proved that the
one for $k[\Delta (d-1,m)]$ is minimal and the corresponding numbers $b_i =\sum_{j} b_{ij}$ are  equal by 
 Equation~(\ref{eqn!keyequalityforbi_from_odd_d_to_d-1}).  This finishes the proof of
Proposition~\ref {prop!keytechnicalforapplication}.
\end {proof}

\subsection  {Examples and implementation}  \label{subs!examplesforKMconstruction}

In this subsection we demonstrate the construction of the cyclic polytope
resolution with a sequence of two examples. First we carry out the
Kustin-Miller complex construction described in
Subsection~\ref{subs!generalities_about_KM_complexes} for the step
passing from the codimension $4$ complete intersection $J_{2,5}$ 
and the Pfaffian $I_{2,5}$ to the codimension $4$ ideal
$I_{2,6}$. In the second step we pass from $J_{4,7}$
and the Pfaffian $I_{4,7}$ to $I_{4,8}$, using that 
$J_{4,7}$ is equal  to $I_{2,6}$ after a change of variables.
At the end of the subsection we give
a link to computer algebra code where we implement our constructions.

Using the notation of Subsection~\ref{subs!generalities_about_KM_complexes},
we will explicitly compute for each
step the auxiliary data $\alpha_{i}$,
$\beta_{i}$, $h_{i}$, $u$ and hence the differentials $f_{i}$
from the input data $a_{i}$ and $b_{i}$.
The ideals $I_{2,5}$ and $I_{4,7}$
are Gorenstein codimension $3$, hence Pfaffian, 
and we will fix below
a certain resolution for each of them.
In addition,  we will also fix  below a certain  Koszul complex resolving 
$J_{2,5} = (z,x_2, \dots  ,x_4)$.

Assume  $q \geq 3$ is an odd integer  and  $M$ is a skew-symmetric
$q \times q $ matrix with entries in a commutative ring. For $1 \leq i \leq  q$,
we denote by $\pf_i M$ the Pfaffian (\cite[Section~3.4]{BH}) 
of the submatrix of $M$ obtained
by deleting the $i$-th row and column of $M$.
The main property of  $\pf_i M$ is that its square is 
the determinant of the corresponding submatrix.

We will use the notation  $R_{(m)} = k[x_1, \dots ,x_m,z]$ introduced
in Section~\ref{sec!combinatorial_interpretation}.
For $d \geq 2$ even,  we denote by $M_d$ the  $(d+3)\times (d+3)$  
skew-symmetric matrix with entries in $R_{(d+3)}$ 
whose $(i,j)$ entry for $i \leq j$ is zero except that 
for $ 1 \leq i \leq d+2$ we have  $(M_d)_{i,i+1} = x_i$ and that
$(M_d)_{1,d+3} = -x_{d+3}$.  It is an easy calculation that 
\[
   I_{d,d+3} =  ( \pf_i (M_d) \bigm| 1 \leq i \leq d+3  ).
\]  
In addition, according to the Buchsbaum-Eisenbud 
theorem  \cite{BE}, 
the minimal graded free resolution of $R_{(d+3)}/I_{d,d+3}$
is given by
\begin{equation}
   0\to R_{(d+3)} \xrightarrow{v_d^t} R_{(d+3)}^{d+3}  
           \xrightarrow{M_d} R_{(d+3)}^{d+3}  \xrightarrow{v_d} R_{(d+3)} 
\label{eqn!resolutionofPfM_d}
\end{equation}
where $v_d $ denote the $1 \times (d+3)$ matrix with $(1,i)$ entry equal
to $(-1)^i \pf_i (M_d)$ and $v_d^t $ denotes the transpose of $v_d$.

We set $R=R_{(5)}$ and fix the following Koszul complex resolution of $R/J_{2,5}$
\begin{equation}
     0\to R \xrightarrow{a_4} R^4   \xrightarrow{a_3} R^6 \xrightarrow{a_2} R^4 
              \xrightarrow{a_1} R 
\label {eqn!resolutionofJ25}             
\end{equation}
where
\small
\[
   a_1  =  \begin{pmatrix}  z & x_3 & x_4 &  x_2  \end{pmatrix},
   \quad \quad  a_2 =  \begin{pmatrix}  x_3 & x_4 & x_2 & 0    & 0    & 0  \\
                           -z  &  0  & 0   & 0    & x_2  & -x_4   \\
                            0  & -z  & 0   & -x_2 &  0   & x_3   \\                      
                            0  &  0  & -z  &  x_4 & -x_3 &  0   
          \end{pmatrix},
\]

\[
a_3 =  \begin{pmatrix}      0  & -x_2 & x_4 & z  & 0  &   0  \\
                            x_2  &  0  & -x_3   & 0  & z  &  0    \\
                            -x_4  &  x_3  & 0   & 0  & 0 & z \\                      
                            0  &  0  & 0  &  x_3   & x_4  & x_2
                            \end{pmatrix}^t,  
      \quad \quad  a_4 =  \begin{pmatrix}  x_3  \\  x_4 \\ x_2  \\ -z     \end{pmatrix}. 
\]

\normalsize
We now discuss the Kustin--Miller complex construction for the
step passing from $(I_{2,5},J_{2,5})$ to $I_{2,6}$,
which corresponds to the  unprojection of
$J_{2,5} \subset R/I_{2,5}$.
We will use as input for the Kustin--Miller complex construction
the resolution  (\ref{eqn!resolutionofJ25}) of $R/J_{2,5}$ and 
the case $d=2$ of (\ref{eqn!resolutionofPfM_d}), which is a 
resolution of  $R/I_{2,5}$. Performing the computations we obtain, 
in the notation of
Subsection~\ref{subs!generalities_about_KM_complexes},  the
complex $C_S$  specified by
$h_{1}=h_{2}=0$, $u=-1$ and the maps%
\begin{align*}
\alpha_{1}  &  :R^{5}\rightarrow R^{4}\text{, }%
{\textstyle\sum\nolimits_{i=1}^{5}}
c_{i}e_{i}\mapsto x_{1}\left(  c_{5}e_{2}+c_{3}e_{3}\right)  +x_{4}c_{1}%
e_{4}+x_{5}\left(  c_{2}e_{2}+c_{4}e_{4}\right) \\
\alpha_{2}  &  :R^{5}\rightarrow R^{6}\text{, }%
{\textstyle\sum\nolimits_{i=1}^{5}}
c_{i}e_{i}\mapsto x_{1}\left(  c_{2}e_{4}+c_{4}e_{6}\right)  +x_{5}c_{3}%
e_{5}\\
\alpha_{3}  &  :R\rightarrow R^{4}\text{, }e_{1}\mapsto x_{1}x_{5}e_{4}%
\end{align*}
and%
\begin{align*}
\beta_{1}  &  :R^{4}\rightarrow R\text{, }%
{\textstyle\sum\nolimits_{i=1}^{4}}
c_{i}e_{i}\mapsto-x_{1}x_{5}c_{1}e_{1}\\
\beta_{2}  &  :R^{6}\rightarrow R^{5}\text{, }%
{\textstyle\sum\nolimits_{i=1}^{6}}
c_{i}e_{i}\mapsto-x_{1}\left(  c_{1}e_{2}+c_{3}e_{4}\right)  -x_{5}c_{2}%
e_{3}\\
\beta_{3}  &  :R^{4}\rightarrow R^{5}\text{, }%
{\textstyle\sum\nolimits_{i=1}^{4}}
c_{i}e_{i}\mapsto-x_{1}\left(  c_{2}e_{3}+c_{1}e_{5}\right)  -x_{4}c_{3}%
e_{1}-x_{5}\left(  c_{1}e_{2}+c_{3}e_{4}\right),
\end{align*}
where  $(e_i)_{1 \leq i \leq q}$ denotes the canonical basis of $R^q$ as $R$-module.
Substituting $x_6$ for $T$  and $0$ for $z$ in the differential 
maps of $C_S$ we get the minimal graded free resolution
of $ R_{(6)}/I_{2,6}$. Moreover, substituting $z$ for $x_1$ in the differential 
maps of the resolution of $R_{(6)}/I_{2,6}$ just constructed we get the minimal
graded free resolution of $R_{(7)}/J_{4,7}$.

We now set $R=R_{(7)}$ and discuss the Kustin--Miller complex 
construction for the step passing from $(I_{4,7},J_{4,7})$ to $I_{4,8}$,
which corresponds to the  unprojection of $J_{4,7} \subset R/I_{4,7}$.
We will use as input for the Kustin--Miller complex construction
the resolution of $R/J_{4,7}$ constructed above and 
the case $d=4$ of (\ref{eqn!resolutionofPfM_d}), which is a 
resolution of  $R/I_{4,7}$. Performing the computations we
obtain, in the notation of
Subsection~\ref{subs!generalities_about_KM_complexes},  the
complex $C_S$  specified by 
$h_{1}%
=h_{2}=0$, $u=-1$ and the maps%
\footnotesize
\begin{align*}
\alpha_{1}  &  :R^{7}\rightarrow R^{9}\text{, }%
{\textstyle\sum\nolimits_{i=1}^{7}}
c_{i}e_{i}\mapsto x_{1}\left(  c_{7}e_{2}+c_{5}e_{7}+c_{3}e_{8}\right)
+x_{6}c_{1}e_{1}+x_{7}\left(  c_{6}e_{1}+c_{2}e_{2}+c_{4}e_{4}\right) \\
\alpha_{2}  &  :R^{7}\rightarrow R^{16}\text{, }%
{\textstyle\sum\nolimits_{i=1}^{7}}
c_{i}e_{i}\mapsto x_{7}\left(  c_{3}e_{3}+c_{5}e_{5}\right)  -x_{1}\left(
c_{2}e_{9}+c_{4}e_{11}-c_{2}e_{12}+c_{6}e_{13}\right) \\
\alpha_{3}  &  :R\rightarrow R^{9}\text{, }e_{1}\mapsto x_{1}x_{7}\left(
x_{5}e_{4}-x_{4}e_{7}-x_{3}e_{9}\right)
\end{align*} 
\normalsize
and%
\footnotesize 
\begin{align*}
\beta_{1}  &  :R^{9}\rightarrow R\text{, }%
{\textstyle\sum\nolimits_{i=1}^{9}}
c_{i}e_{i}\mapsto x_{1}x_{7}\left(  -c_{3}x_{4}-c_{5}x_{3}+c_{6}x_{5}\right)
\\
\beta_{2}  &  :R^{16}\rightarrow R^{7}\text{, }%
{\textstyle\sum\nolimits_{i=1}^{16}}
c_{i}e_{i}\mapsto-x_{1}\left(  c_{1}e_{2}+c_{6}e_{2}+c_{8}e_{4}-c_{2}%
e_{6}\right)  -x_{7}\left(  c_{14}e_{3}+c_{16}e_{5}\right) \\
\beta_{3}  &  :R^{9}\rightarrow R^{7}\text{, }%
{\textstyle\sum\nolimits_{i=1}^{9}}
c_{i}e_{i}\mapsto-x_{6}c_{5}e_{1}-x_{7}\left(  c_{6}e_{2}+c_{8}e_{4}%
+c_{5}e_{6}\right)  +x_{1}\left(  c_{2}e_{3}+c_{1}e_{5}-c_{6}e_{7}\right).
\end{align*} %
\normalsize %
Substituting $x_8$ for $T$  and $0$ for $z$ 
in the differential maps of  $C_S$ we get the minimal graded free 
resolution of $ R_{(8)}/I_{4,8}$.

Under the link \cite{BP2}, 
a  related package for the computer algebra
system Macaulay2 \cite{GS} is available. 
Applying the ideas of the present paper, it
constructs the resolution of the ideal $I_{d,m}$ for  
$d$ even and $m \geq d+1$ 
starting from Koszul complexes  and the
skew-symmetric Buchsbaum-Eisenbud resolution (\ref{eqn!resolutionofPfM_d}) 
of $I_{d,d+3}$.
The functions in the package provide the
user with the option to output all the intermediate data $a_{i}$, $b_{i}$,
$\alpha_{i}$, $\beta_{i}$, $h_{i}$, $u$, $f_{i}$ in addition to the final resolution.

\vspace{0.1 in}
\noindent
\emph{Acknowledgements}. 
The authors wish to thank Christos Athanasiadis for useful discussions and suggestions,
and an anonymous referee for suggestions that
improved the presentation of the material.

\end{document}